\documentclass{article} 

\usepackage[english]{babel} 
\usepackage{amssymb}
\usepackage{amsmath}
\usepackage{txfonts}
\usepackage{mathdots}
\usepackage[classicReIm]{kpfonts}
\usepackage{graphicx}

\begin{document}


\noindent 

\noindent 

\noindent 

\noindent 

\begin{center}
\textbf{{\Large A fast algorithm for the linear programming problem constrained
with the Weighted power mean -- Fuzzy Relational Equalities (WPM-FRE)}}
\end{center}

\noindent 
\begin{center}
\noindent \textbf{Amin Ghodousian${}^{a}$${}^{*}$\footnote{*\ Corresponding\ author$  $Email\ addresses:\ \ $  $a.ghodousian@ut.ac.ir\ (Amin\ Ghodousian).}, Sara Zal${}^{b}$}
\end{center}
\noindent ${}^{a\ }${{\small Faculty of Engineering Science, College of Engineering, University of Tehran, P.O.Box 11365-4563, Tehran, Iran.}}

\noindent ${}^{b}$${}^{\ }${{\small Department of Engineering Science, College of Engineering, University of Tehran, Tehran, Iran.}}

\noindent \textbf{}

\noindent \textbf{Abstract}
\vskip 0.2in
\noindent In this paper, a linear programming problem is investigated in which the feasible region is formed as a special type of fuzzy relational equalities (FRE). In this type of FRE, fuzzy composition is considered as the weighted power mean operator (WPM). Some theoretical properties of the feasible region are derived and some necessary and sufficient conditions are also presented to determine the feasibility of the problem. Moreover, two procedures are proposed for simplifying the problem. Based on some structural properties of the problem, an algorithm is presented to find the optimal solutions and finally, an example is described to illustrate the algorithm. 


\noindent \textbf{}

\noindent \textbf{1. Introduction}

\vskip 0.2in
\noindent In this paper, we study the following linear optimization model whose constraints are formed as a fuzzy system defined by the weighted power mean operator:
\begin{equation} \label{GrindEQ__1_} 
\begin{array}{l} {\min \, \, \, \, c\, x} \\ {\, \, \, \, \, \, \, \, \, \, \, A\, \varphi \, x=b} \\ {\, \, \, \, \, \, \, \, \, \, \, x\in [0,1]^{n} } \end{array} 
\end{equation} 
where $I=\{ \, 1,2,...,m\, \} $ ,  $J=\{ \, 1,2,...,n\, \} $, $A=(a_{ij} )_{m\, \times n} $, $0\le a_{ij} \le 1$ ($\forall i\in I$\textbf{ }and $\forall j\in J$), is a fuzzy matrix, $b=(b_{i}^{} )_{m\, \times 1} $,$0\le b_{i} \le 1$ ($\forall i\in I$), is an $m$--dimensional fuzzy vector, and ``$\varphi $'' is the max-weighted power mean composition, that is, $\varphi \, (x,y)=\left(w\, x^{\, p} +(1-w)y^{\, p} \right)^{{1\, \mathord{\left/ {\vphantom {1\,  p}} \right. \kern-\nulldelimiterspace} p} } $.

\noindent Furthermore, let \textbf{$S(A,b)$ }denote the feasible solutions sets of problem (1), that is, \textbf{$S(A,b)=\left\{x\in [0,1]^{n} \, \, :\, \, \, A\, \varphi \, x=b\right\}$}. Additionally, if\textbf{ }$a_{i} $\textbf{ }denotes the \textbf{$i$}'th row of matrix $A$, then problem(1) can be also expressed as follows:
\begin{equation} \label{GrindEQ__2_} 
\begin{array}{l} {\min \, \, \, \, cx} \\ {\, \, \, \, \, \, \, \, \, \, \, \varphi \, (a_{i} \, ,x)=b_{i}^{} \, \, \, ,\, \, \, i\in I} \\ {\, \, \, \, \, \, \, \, \, \, \, x\in [0,1]^{n} } \end{array} 
\end{equation} 
where the constraints mean $\varphi \, (a_{i} ,x)={\mathop{\max }\limits_{j\in J}} \{ \, \varphi \, (a_{ij} ,x_{j} )\, \} =b_{i} $ ($\forall i\in I$) and \textbf{ }$\varphi \, (a_{ij} ,x_{j} )=\left(w\, a_{ij}^{\, p} +(1-w)x_{j}^{\, p} \right)^{{1\, \mathord{\left/ {\vphantom {1\,  p}} \right. \kern-\nulldelimiterspace} p} } $.\textbf{                                  }                                                                                                                    

\noindent The theory of fuzzy relational equations (FRE) was firstly proposed by Sanchez and applied in problems of the medical diagnosis [39]. Nowadays, it is well known that many issues associated with a body knowledge can be treated as FRE problems [35]. In addition to the preceding applications, FRE theory has been applied in many fields, including fuzzy control, discrete dynamic systems, prediction of fuzzy systems, fuzzy decision making, fuzzy pattern recognition, fuzzy clustering, image compression and reconstruction, fuzzy information retrieval, and so on. Generally, when inference rules and their consequences are known, the problem of determining antecedents is reduced to solving an FRE [25,33]. 

\noindent The solvability determination and the finding of solutions set are the primary (and the most fundamental) subject concerning with FRE problems. Actually, The solution set of FRE is often a non-convex set that is completely determined by one maximum solution and a finite number of minimal solutions [5]. This non-convexity property is one of two bottlenecks making major contribution to the increase of complexity in problems that are related to FRE, especially in the optimization problems subjected to a system of fuzzy relations. The other bottleneck is concerned with detecting the minimal solutions for FREs [2]. Markovskii showed that solving max-product FRE is closely related to the covering problem which is an NP-hard problem [32].  In fact, the same result holds true for a more general t-norms instead of the minimum and product operators [2,3,12,13,15,16,28,29,32]. 

\noindent Over the last decades, the solvability of FRE defined with different max-t compositions have been investigated by many researchers [15,16,34,36,37,40,42,43,\\*45,48,51]. Moreover, some researchers introduced and improved theoretical aspects and applications of fuzzy relational inequalities (FRI) [12 -- 14,17,18,26\\*,50]. Li and Yang [26] studied a FRI with addition-min composition and presented an algorithm to search for minimal solutions. Ghodousian et al. [13] focused on the algebraic structure of two fuzzy relational inequalities \textbf{$A\varphi x\le b^{1} $}and \textbf{$D\varphi x\ge b^{2} $}, and studied a mixed fuzzy system formed by the two preceding FRIs, where $\varphi $is an operator with (closed) convex solutions. 

\noindent The problem of optimization subject to FRE and FRI is one of the most interesting and on-going research topic among the problems related to FRE and FRI theory [1,8,11 -- 16,23,27,30,38,41,46,50]. Fang and Li [9] converted a linear optimization problem subjected to FRE constraints with max-min operation into an integer programming problem and solved it by branch and bound method using jump-tracking technique. In [23] an application of optimizing the linear objective with max-min composition was employed for the streaming media provider. Wu et al. [44] improved the method used by Fang and Li, by decreasing the search domain. The topic of the linear optimization problem was also investigated with max-product operation [20,31]. Loetamonphong and Fang defined two sub-problems by separating negative and non-negative coefficients in the objective function and then obtained the optimal solution by combining those of the two sub-problems [31]. Also, in [20] some necessary conditions of the feasibility and simplification techniques were presented for solving FRE with max-product composition. Moreover, some generalizations of the linear optimization with respect to FRE have been studied with the replacement of max-min and max-product compositions with different fuzzy compositions such as max-average composition [46] and max-t-norm composition [15,16,21,27,41]. 

\noindent Recently, many interesting generalizations of the linear programming subject to a system of fuzzy relations have been introduced and developed based on composite operations used in FRE, fuzzy relations used in the definition of the constraints, some developments on the objective function of the problems and other ideas [6,10,15,16,18,24,30,47].  For example, Dempe and Ruziyeva [4] generalized the fuzzy linear optimization problem by considering fuzzy coefficients. 

\noindent The optimization problem subjected to various versions of FRI could be found in the literature as well [12 -- 14,17,18,49,50]. Xiao et al. [50] introduced the latticized linear programming problem subject to max-product fuzzy relation inequalities. Ghodousian et al. [12] introduced a system of fuzzy relational inequalities with fuzzy constraints (FRI-FC) in which the constraints were defined with max-min composition. 

\noindent In this paper, an algorithm is proposed to find all the optimal solutions of problem (1). Firstly, we describe some structural details of WPM-FREs such as the theoretical properties of WPM-fuzzy equalities and necessary and sufficient conditions for the feasibility of the problem. Then, the feasible region is completely determined by a finite number of convex cells. Additionally, some simplification processes are introduced to reduce the problem. Finally, an algorithm is presented to solve the main problem.  

\noindent The remainder of the paper is organized as follows. Section 2 gives some basic results on the WPM-fuzzy equalities. Also, some feasibility conditions are derived. In section 3, the feasible region is characterized in terms of a finite number of closed convex cells. In section 4, some simplification rules are presented. These rules convert the problem into an equivalent one that is easier to solve. The optimal solution of the problem is described in Section 5 and, finally in section 6 an example is presented to illustrate the algorithm.

\noindent 
\vskip 0.2in
\noindent \textbf{2. Basic properties of WPM -- FREs }
\vskip 0.2in
\noindent In this section, the structural properties of each fuzzy equation $\varphi \, (a_{i} ,x)=b_{i} $ is investigated and its solutions are found. Let $S(a_{i} ,b_{i} )$\textbf{ }denote the feasible solutions set of \textbf{$i$}`th equation, that is, $S(a_{i} ,b_{i} )=\left\{x\in [0,1]^{n} \, \, :\, \, \, \varphi \, (a_{i} ,x)=b_{i} \right\}$. So, $S(A,b)=\bigcap _{i\in I}S(a_{i} ,b_{i} ) $. 

\noindent \textbf{}

\noindent \textbf{Lemma 1.} Let $i\in I$, $j_{0} \in J$ and $a_{ij_{0} } >\frac{b_{i} }{\sqrt[{p}]{w} } $. Then, $S(a_{i} ,b_{i} )=\emptyset $.
\vskip 0.2in
\noindent \textbf{Proof.} Since $\varphi $ is an increasing function on $[0,1]^{2} $ in both variables, we note that $\varphi \, (a_{ij_{0} } ,x_{j} )>\varphi \, (\, {b_{i} \mathord{\left/ {\vphantom {b_{i}  \sqrt[{p}]{w} }} \right. \kern-\nulldelimiterspace} \sqrt[{p}]{w} } \, ,\, x_{j} )=\left(b_{i}^{\, p} +(1-w)x_{j_{0} }^{\, p} \right)^{{1\, \mathord{\left/ {\vphantom {1\,  p}} \right. \kern-\nulldelimiterspace} p} } \ge b_{i} $. Thus, for each $x\in [0,1]^{n} $ we have \textbf{$\varphi \, (a_{i} ,x)={\mathop{\max }\limits_{j\in J}} \{ \, \varphi \, (a_{ij} ,x_{j} )\, \} \ge \varphi \, (a_{ij_{0} } ,x_{j_{0} } )>b_{i} $}. Hence, $x\notin S(a_{i} ,b_{i} )$, $\forall x\in [0,1]^{n} $.  $\square$

\vskip 0.2in
\noindent 

\noindent \textbf{Lemma 2. }Let $a_{ij_{0} } \le \frac{b_{i} }{\sqrt[{p}]{w} } $ for some $i\in I$ and $j_{0} \in J$. If $b_{i}^{\, p} \ge 1-w$ and $a_{ij_{0} } <\left(\frac{b_{i}^{\, p} +w-1}{w} \right)^{{1\, \mathord{\left/ {\vphantom {1\,  p}} \right. \kern-\nulldelimiterspace} p} } $, then $\varphi \, (a_{ij_{0} } ,x_{j_{0} } )<b_{i} $, $\forall x_{j_{0} } \in [0,1]$. 
\vskip 0.2in
\noindent \textbf{Proof.} Since  $b_{i}^{\, p} \ge 1-w$, then $[(b_{i}^{\, p} +w-1)/w]^{{1\, \mathord{\left/ {\vphantom {1\,  p}} \right. \kern-\nulldelimiterspace} p} } \ge 0$.  Now, the result follows from the relations $\varphi \, (a_{ij_{0} } ,x_{j_{0} } )<\varphi \, (\, \, [(b_{i}^{\, p} +w-1)/w]^{\, {1\, \mathord{\left/ {\vphantom {1\,  p}} \right. \kern-\nulldelimiterspace} p} } ,1\, )=b_{i} $. $\square$ 

\noindent 
\vskip 0.2in
\noindent \textbf{Lemma 3. }Let $a_{ij_{0} } \le \frac{b_{i} }{\sqrt[{p}]{w} } $ for some $i\in I$ and $j_{0} \in J$. Also, suppose that either $b_{i}^{\, p} <1-w$ or $a_{ij_{0} } \ge \left(\frac{b_{i}^{\, p} +w-1}{w} \right)^{{1\, \mathord{\left/ {\vphantom {1\,  p}} \right. \kern-\nulldelimiterspace} p} } $. Then, $x_{j_{0} } =\left(\frac{b_{i}^{\, p} -w\, a_{ij_{0} }^{\, p} }{1-w} \right)^{{1\, \mathord{\left/ {\vphantom {1\,  p}} \right. \kern-\nulldelimiterspace} p} } $ is the unique solution to the equality $\varphi \, (a_{ij_{0} } ,x_{j_{0} } )=b_{i} $. 
\vskip 0.2in
\noindent \textbf{Proof.} It is easy to verify that $\varphi \, (a_{ij_{0} } ,x_{j_{0} } )=b_{i} $. Now, since $\varphi $ is an increasing function, we have $\varphi \, (a_{ij_{0} } ,x_{j} )>b_{i} $ if $x_{j} >[(b_{i}^{\, p} -w\, a_{ij_{0} }^{\, p} )/(1-w)]^{{1\, \mathord{\left/ {\vphantom {1\,  p}} \right. \kern-\nulldelimiterspace} p} } $ and $\varphi \, (a_{ij_{0} } ,x_{j} )<b_{i} $ if $x_{j} <[(b_{i}^{\, p} -w\, a_{ij_{0} }^{\, p} )/(1-w)]^{{1\, \mathord{\left/ {\vphantom {1\,  p}} \right. \kern-\nulldelimiterspace} p} } $. $\square$

\noindent \textbf{}

\noindent From Lemmas 1, 2 and 3, the following theorem is resulted that gives a necessary and sufficient condition for the feasibility of the set $S(a_{i} ,b_{i} )$.

\noindent 
\vskip 0.2in
\noindent \textbf{Theorem 1. }For a fixed $i\in I$, $S(a_{i} ,b_{i} )\ne \emptyset $\textbf{ }if and only if 

\noindent (a) $a_{ij} \le {b_{i} \mathord{\left/ {\vphantom {b_{i}  \sqrt[{p}]{w} }} \right. \kern-\nulldelimiterspace} \sqrt[{p}]{w} } $, $\forall j\in J$.

\noindent (b) There exist some $j_{0} \in J$ such that either $b_{i}^{\, p} <1-w$ or $a_{ij_{0} } \ge [(b_{i}^{\, p} +w-1)/(w)]^{{\, 1\, \mathord{\left/ {\vphantom {\, 1\,  p}} \right. \kern-\nulldelimiterspace} p} } $.
\vskip 0.2in
\noindent \textbf{Definition 1. }For an arbitrary fixed $i\in I$, let $J^{\, -} (i)=\left\{j\in J\, \, :\, \, \, a_{ij} >{b_{i} \mathord{\left/ {\vphantom {b_{i}  \sqrt[{p}]{w} }} \right. \kern-\nulldelimiterspace} \sqrt[{p}]{w} } \right\}$. Additionally, define  $J^{\, \infty } (i)=\left\{j\in J\, \, :\, \, \, b_{i}^{\, p} \ge 1-w\, \, ,\, \, a_{ij} <[(b_{i}^{\, p} +w-1)/w]^{{1\, \mathord{\left/ {\vphantom {1\,  p}} \right. \kern-\nulldelimiterspace} p} } \right\}$ and $J(i)=J-\left\{J^{\, -} (i)\bigcup J^{\, \infty } (i)\right\}$. 

\noindent 

\noindent According to Theorem 1, the following corollary is directly attained. This corollary characterizes all the feasible solutions of $S(a_{i} ,b_{i}^{} )$.

\noindent 
\vskip 0.2in
\noindent \textbf{Corollary 1. }$x'\in S(a_{i} ,b_{i}^{} )$if and only if $J^{\, -} (i)=\emptyset $ , $J(i)\ne \emptyset $ and 

\noindent (a) $x'_{j} \in [0,1]$, $\forall j\in J^{\, \infty } (i)$. 

\noindent (b) $x'_{j} \le [(b_{i}^{\, p} -w\, a_{ij}^{\, p} )/(1-w)]^{{1\, \mathord{\left/ {\vphantom {1\,  p}} \right. \kern-\nulldelimiterspace} p} } $, $\forall j\in J(i)$.

\noindent (c) There exist at least some $j_{0} \in J(i)$ such that $x'_{j_{0} } =[(b_{i}^{\, p} -w\, a_{ij_{0} }^{\, p} )/(1-w)]^{{1\, \mathord{\left/ {\vphantom {1\,  p}} \right. \kern-\nulldelimiterspace} p} } $.

\noindent 
\vskip 0.2in
\noindent \textbf{Definition 2. }Suppose that $S(a_{i} ,b_{i}^{} )\ne \emptyset $(hence, $J^{\, -} (i)=\emptyset $ from Corollary 1). Define $\overline{X}(i)\in [0,1]^{n} $ such that 
\[\overline{X}(i)_{j} =\left\{\begin{array}{l} {\, \left(\frac{b_{i}^{\, p} -w\, a_{ij}^{\, p} }{1-w} \right)^{{1\, \mathord{\left/ {\vphantom {1\,  p}} \right. \kern-\nulldelimiterspace} p} } \, \, \, \, \, \, \, ,\, if\, \, j\in J(i)} \\ {\, \, \, \, \, \, \, \, \, \, \, \, \, 1\, \, \, \, \, \, \, \, \, \, \, \, \, \, \, \, \, \, \, \, \, \, \, \, \, \, ,\, if\, \, j\in J^{\, \infty } (i)} \end{array}\right. \] 
\textbf{}

\noindent \textbf{Theorem 2. }Suppose that $S(a_{i} ,b_{i}^{} )\ne \emptyset $. Then, $\overline{X}(i)$ is the maximum solution of $S(a_{i} ,b_{i}^{} )$. 
\vskip 0.2in
\noindent \textbf{Proof. }Since $S(a_{i} ,b_{i}^{} )\ne \emptyset $, then $J^{\, -} (i)=\emptyset $. Based on Corollary 1, $\overline{X}(i)\in S(a_{i} ,b_{i}^{} )$. Suppose that $x'\in S(a_{i} ,b_{i}^{} )$. So, from Corollary 1, $x'_{j} \le [(b_{i}^{\, p} -w\, a_{ij}^{\, p} )/(1-w)]^{{1\, \mathord{\left/ {\vphantom {1\,  p}} \right. \kern-\nulldelimiterspace} p} } $, $\forall j\in J(i)$, and $x'_{j} \le 1$, $\forall j\in J^{\, \infty } (i)$. Therefore, $x'_{j} \le \overline{X}(i)_{j} $, $\forall j\in J$. $\square$

\noindent 
\vskip 0.2in
\noindent \textbf{Definition 3. }Let $i\in I$ and $S(a_{i} ,b_{i}^{} )\ne \emptyset $. For each $j\in J(i)$, define $\underline{X}(i,j)\in [0,1]^{n} $ such that
\[\underline{X}(i,j)_{k} =\left\{\begin{array}{l} {\, \left(\frac{b_{i}^{\, p} -w\, a_{ij}^{\, p} }{1-w} \right)^{{1\, \mathord{\left/ {\vphantom {1\,  p}} \right. \kern-\nulldelimiterspace} p} } \, \, \, \, \, \, \, ,\, k=j} \\ {\, \, \, \, \, \, \, \, \, \, \, \, \, 0\, \, \, \, \, \, \, \, \, \, \, \, \, \, \, \, \, \, \, \, \, \, \, \, \, \, ,\, otherwise} \end{array}\right. \]

\noindent \textbf{Remark 1.} Suppose that $S(a_{i} ,b_{i}^{} )\ne \emptyset $ and $j\in J(i)$. Then, from Definitions 2 and 3, we have $\overline{X}(i)_{j} =\underline{X}(i,j)_{j} $.

\noindent 
\vskip 0.2in
\noindent \textbf{Theorem 3. }Suppose that $S(a_{i} ,b_{i}^{} )\ne \emptyset $ and $j_{0} \in J(i)$. Then, $\underline{X}(i,j_{0} )$ is a minimal solution of $S(a_{i} ,b_{i}^{} )$. 
\vskip 0.2in
\noindent \textbf{Proof. }From Corollary 1, $\underline{X}(i,j_{0} )\in S(a_{i} ,b_{i}^{} )$. Suppose that $x'\in S(a_{i} ,b_{i}^{} )$ , $x'\le \underline{X}(i,j_{0} )$and $x'\ne \underline{X}(i,j_{0} )$. So, $x'_{j} \le \underline{X}(i,j_{0} )_{j} $, $\forall j\in J$ and $x'\ne \underline{X}(i,j_{0} )$. Therefore, $x'_{j} =0$, $\forall j\in J-\{ j_{0} \} $, and $x'_{j_{0} } <[(b_{i}^{\, p} -w\, a_{ij_{0} }^{\, p} )/(1-w)]^{{1\, \mathord{\left/ {\vphantom {1\,  p}} \right. \kern-\nulldelimiterspace} p} } $. Hence, from Lemmas 1, 2 and 3 we have \textbf{$\varphi \, (a_{i} ,x')=\max \, \left\{{\mathop{\max }\limits_{j\in J-\{ j_{0} \} }} \{ \, \varphi \, (a_{ij} ,x'_{j} )\, \} \, ,\, \varphi \, (a_{ij_{0} } ,x'_{j_{0} } )\right\}=\varphi \, (a_{ij_{0} } ,x'_{j_{0} } )<b_{i} $} that contradicts $x'\in S(a_{i} ,b_{i}^{} )$. $\square$

\noindent 

\noindent The following theorem shows that $S(a_{i} ,b_{i} )$ can be stated in terms of the unique maximum solution and a finite number of minimal solutions. 

\noindent 
\vskip 0.2in
\noindent \textbf{Theorem 4. }$S(a_{i} ,b_{i} )=\bigcup _{j\in J(i)}[\underline{X}(i,j)\, ,\, \overline{X}(i)\, ] $. 
\vskip 0.2in
\noindent \textbf{Proof. }Let $x'\in S(a_{i} ,b_{i}^{} )$. From Theorem 2, $x'\le \overline{X}(i)$. Furthermore, there exist at least some $j_{0} \in J(i)$ such that $x'_{j_{0} } =[(b_{i}^{\, p} -w\, a_{ij_{0} }^{\, p} )/(1-w)]^{{1\, \mathord{\left/ {\vphantom {1\,  p}} \right. \kern-\nulldelimiterspace} p} } $ (Corollary 1). Thus, from Definition 3 we have $\underline{X}(i,j_{0} )\le x'$. Hence, $x'\in [\underline{X}(i,j_{0} )\, ,\, \overline{X}(i)\, ]$. Conversely, let $x'\in \bigcup _{j\in J(i)}[\underline{X}(i,j)\, ,\, \overline{X}(i)\, ] $. Therefore, \textbf{$\varphi \, (a_{ij} ,x'_{j} )\le \varphi \, (a_{ij} ,\overline{X}(i)_{j} )\le b_{i} $}, $\forall j\in J$. Moreover, there exists some $j_{0} \in J(i)$ such that $x'\in [\underline{X}(i,j_{0} )\, ,\, \overline{X}(i)\, ]$. So, Remark 1 implies $x'_{j_{0} } =\underline{X}(i,j_{0} )_{j_{0} } =\overline{X}(i)_{j_{0} } $ and therefore, \textbf{$\varphi \, (a_{ij_{0} } ,x'_{j_{0} } )=b_{i} $}. Thus, we have
\[\varphi \, (a_{i} ,x')={\mathop{\max }\limits_{j\in J}} \{ \, \varphi \, (a_{ij} ,x'_{j} )\, \} =\max \, \left\{{\mathop{\max }\limits_{j\in J-\{ j_{0} \} }} \{ \, \varphi \, (a_{ij} ,x'_{j} )\, \} \, ,\, \varphi \, (a_{ij_{0} } ,x'_{j_{0} } )\right\}=\varphi \, (a_{ij_{0} } ,x'_{j_{0} } )=b_{i} \] 
which implies that $x'\in S(a_{i} ,b_{i}^{} )$. $\square$

\noindent 
\vskip 0.2in
\noindent \textbf{3. Feasible region of Problem (1) }
\vskip 0.2in
\noindent In this section, a necessary and sufficient condition is derived to determine the feasibility of the main problem. 

\noindent 
\vskip 0.2in
\noindent \textbf{Definition 4. }Let $\overline{X}(i)$ be as in Definition 2, \textbf{$\forall i\in I$}. We define $\overline{X}={\mathop{\min }\limits_{i\in I}} \, \left\{\overline{X}(i)\right\}$. 

\noindent \textbf{}

\noindent \textbf{Definition 5. }Let $e:I\to \bigcup _{i\in I}J(i) $\textbf{ }so that $e(i)\in J(i)$, $\forall i\in I$, and let $E$ be the set of all vectors $e$. For the sake of convenience, we represent each $e\in E$ as an $m$--dimensional vector $e=[j_{1} ,j_{2} ,...,j_{m} ]$ in which $j_{k} =e(k)$, $k=1,2,...,m$.

\noindent \textbf{}

\noindent \textbf{Definition 6. }Let $e=[j_{1} ,j_{2} ,...,j_{m} ]\in E$. We define $\underline{X}(e)\in [0,1]^{n} $ such that $\underline{X}(e)_{j} ={\mathop{\max }\limits_{i\in I}} \left\{\underline{X}(i,e(i))_{j} \right\}={\mathop{\max }\limits_{i\in I}} \left\{\underline{X}(i,j_{i} )_{j} \right\}$, $\forall j\in J$.\textbf{  }

\noindent \textbf{}

\noindent The following theorem indicates that the feasible region of problem 1 is completely found by a finite number of closed convex cells. 

\noindent \textbf{}

\noindent \textbf{Theorem 5. $S(A,b)=\bigcup _{e\in E}[\underline{X}(e),\overline{X}] $}.  
\vskip 0.2in
\noindent \textbf{Proof. }Since $S(A,b)=\bigcap _{i\in I}S(a_{i} ,b_{i}^{} ) $, from Theorem 4 we have $S(A,b)=\bigcap _{i\in I}\bigcup _{j\in J_{i} }[\underline{X}(i,j)\, ,\, \overline{X}(i)\, ]  $. So, $S(A,b)=\bigcup _{e\in E}\bigcap _{i\in I}[\underline{X}(i,e(i))\, ,\, \overline{X}(i)\, ]  $ (see Definitions 5 and 6), i.e., $S(A,b)=\bigcup _{e\in E}\, \left[{\mathop{\max }\limits_{i\in I}} \, \left\{\underline{X}(i,e(i))\right\}\, ,\, {\mathop{\min }\limits_{i\in I}} \, \left\{\overline{X}(i)\right\}\right] $. Now, the result follows from Definitions 4 and 6. $\square$

\noindent \textbf{}

\noindent The following Corollary gives a simple necessary and sufficient condition for the feasibility of \textbf{$S(A,b)$}. \textbf{}

\noindent \textbf{}

\noindent \textbf{Corollary 2.} $S(A,b)\ne \emptyset $ iff $\overline{X}\in S(A,b)$. 

\noindent \textbf{}
\newpage
\noindent \textbf{4. Simplification techniques}
\vskip 0.2in
\noindent In practice, there are often some components of matrix \textbf{$A$ }that have no effect on the solutions to problem (1). Therefore, we can simplify the problem by changing the values of these components to zeros. For this reason, various simplification processes have been proposed by researchers. We refer the interesting reader to [13] where a brief review of such these processes is given. Here, we present two simplification techniques based on the weighted power mean operator.

\noindent 
\vskip 0.2in
\noindent \textbf{Definition 7. }If a value changing in an element, say $a_{ij} $, of a given fuzzy relation matrix $A$ has no effect on the solutions of problem (1), this value changing is said to be an equivalence operation.

\noindent 
\vskip 0.2in
\noindent \textbf{Corollary 3. }Suppose that $\varphi \, (a_{ij_{0} } ,x_{j_{0} } )<b_{i} $,$\forall x\in S(A,b)$. In this case, it is obvious that $\mathop{\max }\limits_{j=1}^{n} \, \left\{\varphi \, (a_{ij} ,x_{j} )\right\}=b_{i} $\textbf{ }is equivalent to $\mathop{\max }\limits_{\begin{array}{l} {j=1} \\ {j\ne j_{0} } \end{array}}^{n} \, \left\{\varphi \, (a_{ij} ,x_{j} )\right\}=b_{i} $,\textbf{ }that is, ``resetting $a_{ij_{0} } $ to zero'' has no effect on the solutions of problem (1) (since component $a_{ij_{0} } $ only appears in the \textbf{$i$}`th constraint of problem (1). Therefore, if $\varphi \, (a_{ij_{0} } ,x_{j_{0} } )<b_{i} $,$\forall x\in S(A,b)$, then ``resetting $a_{ij_{0} } $ to zero'' is an equivalence operation. 

\noindent \textbf{}

\noindent \textbf{Lemma 4 (first simplification). }Suppose that $j_{0} \in J^{\, \infty } (i)$, for some $i\in I$ and $j_{0} \in J$. Then, ``resetting $a_{ij_{0} } $ to zero'' is an equivalence operation.

\noindent \textbf{Proof.} The proof is directly resulted from Lemma 2. $\square$\textbf{}

\noindent \textbf{}

\noindent \textbf{Lemma 5 (second simplification). }Suppose that $j_{0} \in J(k)$, where $k\in I$ and $j_{0} \in J$. If there exists some $i\in I$ ($i\ne k$) such that $j_{0} \in J(i)$ and $b_{i}^{\, p} -b_{k}^{\, p} <w\, (a_{ij_{0} }^{\, p} -a_{k\, j_{0} }^{\, p} )$, then ``resetting $a_{k\, j_{0} } $ to zero'' is an equivalence operation.
\vskip 0.2in
\noindent  \textbf{Proof. }We show that $\varphi \, (a_{k\, j_{0} } ,x_{j_{0} } )<b_{k} $,$\forall x\in S(A,b)$. Consider an arbitrary feasible solution $x\in S(A,b)$. Since $x\in S(A,b)$, it turns out that $\varphi \, (a_{k\, j_{0} } ,x_{j_{0} } )>b_{k} $\textbf{ }never holds. So, assume that $\varphi \, (a_{k\, j_{0} } ,x_{j_{0} } )=b_{k} $. Since $j_{0} \in J(k)$, from lemma 3 we conclude that $x_{j_{0} } =[(b_{k}^{\, p} -w\, a_{k\, j_{0} }^{\, p} )/(1-w)]^{{1\, \mathord{\left/ {\vphantom {1\,  p}} \right. \kern-\nulldelimiterspace} p} } $. On the other hand, inequality $b_{i}^{\, p} -b_{k}^{\, p} <w\, (a_{ij_{0} }^{\, p} -a_{k\, j_{0} }^{\, p} )$ implies that $[(b_{i}^{\, p} -w\, a_{ij_{0} }^{\, p} )/(1-w)]^{{1\, \mathord{\left/ {\vphantom {1\,  p}} \right. \kern-\nulldelimiterspace} p} } <[(b_{k}^{\, p} -w\, a_{k\, j_{0} }^{\, p} )/(1-w)]^{{1\, \mathord{\left/ {\vphantom {1\,  p}} \right. \kern-\nulldelimiterspace} p} } $. So, according to Definitions 3 and 4, $\overline{X}_{j_{0} } \le [(b_{i}^{\, p} -w\, a_{ij_{0} }^{\, p} )/(1-w)]^{{1\, \mathord{\left/ {\vphantom {1\,  p}} \right. \kern-\nulldelimiterspace} p} } <x_{j_{0} } $. Therefore, \textbf{$x\notin \bigcup _{e\in E}[\underline{X}(e),\overline{X}] $} that means \textbf{$x\notin S(A,b)$} (Theorem 5). $\square$

\noindent \textbf{}

\noindent \textbf{5. Resolution of Problem (1) }
\vskip 0.2in
\noindent It can be easily verified that \textbf{$\overline{X}$} is the optimal solution for \textbf{$\min \, \, \, \, \left\{Z_{1} =\sum _{j=1}^{n}c_{j}^{-} x_{j}  \, \, :\, \, A\, \varphi \, x=b\, \, ,\, \, x\in [0,1]^{n} \right\}$}, and the optimal solution for \textbf{$\min \, \, \, \, \left\{Z_{2} =\sum _{j=1}^{n}c_{j}^{+} x_{j}  \, \, :\, \, A\, \varphi \, x=b\, \, ,\, \, x\in [0,1]^{n} \right\}$} is \textbf{$\underline{X}(e^{*} )$} for some \textbf{$e^{*} \in E$}, where \textbf{$c_{j}^{+} =\max \{ c_{j} ,0\} $ }and \textbf{$c_{j}^{-} =\min \{ c_{j} ,0\} $ }for \textbf{$j=1,2,...,n$ }[9,13,19,28]. According to the foregoing results, the following theorem shows that the optimal solution of Problem (1) can be obtained by the combination of \textbf{$\overline{X}$} and \textbf{$\underline{X}(e^{*} )$}.

\noindent 
\vskip 0.2in
\noindent \textbf{Theorem 6. }Suppose that $S(A,b)\ne \emptyset $, and $\overline{X}$\textbf{ }and $\underline{X}(e^{*} )$ are the optimal solutions of sub-problems $Z_{1} $ and $Z_{2} $, respectively. Then $c^{T} x^{*} $\textbf{ }is the lower bound of the optimal objective function in (1), where $x^{*} =[x_{1}^{*} ,x_{2}^{*} ,...,x_{n}^{*} ]$\textbf{ }is defined as follows:
\[x_{j}^{*} =\left\{\begin{array}{cc} {\overline{X}_{j} } & {c_{j} <0} \\ {\underline{X}(e^{*} )_{j} } & {c_{j} \ge 0} \end{array}\right. \] 
for $j=1,2,...,n$.
\vskip 0.2in
\noindent \textbf{Proof. }For a general case, see the proof of Theorem 4.1 in [13]. $\square$

\noindent \textbf{}

\noindent \textbf{Corollary 4. }Suppose that $S(A,b)\ne \emptyset $. Then, $x^{*} $\textbf{ }as defined in Theorem 5, is the optimal solution of problem (1).
\vskip 0.2in
\noindent \textbf{Proof.} According to the definition of vector $x^{*} $, we have $\underline{X}(e^{*} )_{j} \le x_{j}^{*} \le \overline{X}_{j} $, $\forall j\in J$, which implies $x^{*} \in \bigcup _{e\in E}[\underline{X}(e),\overline{X}] =S(A,b)$. $\square$

\noindent \textbf{}

\noindent \textbf{6. Numerical example}
\vskip 0.2in
\noindent Consider the following linear programming problem constrained with a fuzzy system defined by the weighted power mean operator: 

\noindent \textbf{}
\[\begin{array}{l} {\min \, \, \, \, Z=-{\rm 7.6582}\, x_{1} -{\rm 2.029}\, x_{2} +{\rm 6.6277}\, x_{3} -{\rm 6.3}\, x_{4} {\rm +0.0157}\, x_{5} -{\rm 7.4737}\, x_{6} {\rm +7.2926}\, x_{7} } \\ {} \\ {\, \, \, \, \, \, \, \, \, \, \, \left[\begin{array}{l} {{\rm 0.6763\; \; \; 0.8969\; \; \; 0.8403\; \; \; 0.3000\; \; \; 0.0710\; \; \; 0.0758\; \; \; 0.3529}} \\ {{\rm 0.3362\; \; \; 0.2721\; \; \; 0.1956\; \; \; 0.3396\; \; \; 0.0101\; \; \; 0.2557\; \; \; 0.1193}} \\ {{\rm 0.1637\; \; \; 0.5426\; \; \; 0.2534\; \; \; 0.3701\; \; \; 0.4916\; \; \; 0.5761\; \; \; 0.2454}} \\ {{\rm 0.5161\; \; \; 0.1330\; \; \; 0.9090\; \; \; 0.1477\; \; \; 0.3827\; \; \; 0.7212\; \; \; 0.2452}} \\ {{\rm 0.2319\; \; \; 0.8371\; \; \; 0.1275\; \; \; 0.8609\; \; \; 0.5201\; \; \; 0.6163\; \; \; 0.0654}} \end{array}\right]\varphi x=\left[\begin{array}{l} {{\rm 0.8657}} \\ {{\rm 0.6520}} \\ {{\rm 0.6926}} \\ {{\rm 0.8833}} \\ {{\rm 0.8350}} \end{array}\right]} \\ {\, \, \, \, \, \, \, \, \, \, \, x\in [0,1]^{7} } \end{array}\] 
where $\left|\, I\, \right|=5$, $\left|\, J\, \right|=7$ and $\varphi \, (x,y)=\left(w\, x^{\, p} +(1-w)y^{\, p} \right)^{{1\, \mathord{\left/ {\vphantom {1\,  p}} \right. \kern-\nulldelimiterspace} p} } $ in which $w=3/4$ and $p=3$.  Moreover, \textbf{$Z_{1} =-{\rm 7.6582}\, x_{1} -{\rm 2.029}\, x_{2} -{\rm 6.3}\, x_{4} -{\rm 7.4737}\, x_{6} $} and \textbf{$Z_{2} ={\rm 6.6277}\, x_{3} {\rm +0.0157}\, x_{5} {\rm +7.2926}\, x_{7} $}. For each $i\in I$, we have $J^{\, -} (i)=\emptyset $. Also, $J(1)=\{ 2,3\} $, $J\eqref{GrindEQ__2_}=\{ 1,4\} $, $J(3)=\{ 2,5,6\} $, $J(4)=\{ 3\} $ and $J(5)=\{ 2,4\} $. Therefore, by Theorem 1, $S(a_{i} ,b_{i} )\ne \emptyset $, $\forall i\in I$. According to Definition 2, the maximum solutions of $S(a_{i} ,b_{i} )\ne \emptyset $, $\forall i\in I$, are attained as follows:
\[\begin{array}{l} {\overline{X}(1)=[{\rm 1\; ,\; 0.7552\; ,\; 0.9341\; ,\; 1\; ,\; 1\; ,\; 1\; ,\; 1}]} \\ {\overline{X}(2)=[{\rm 0.9982\; ,\; 1\; ,\; 1\; ,\; 0.9970\; ,\; 1\; ,\; 1\; ,\; 1}]} \\ {\overline{X}(3)=[{\rm 1\; ,\; 0.9471\; ,\; 1\; ,\; 1\; ,\; 0.9908\; ,\; 0.9107\; ,\; 1}]} \\ {\overline{X}(4)=[{\rm 1\; ,\; 1\; ,\; 0.7955\; ,\; 1\; ,\; 1\; ,\; 1\; ,\; 1}]} \\ {\overline{X}(5)=[{\rm 1\; ,\; 0.8286\; ,\; 1\; ,\; 0.7456\; ,\; 1\; ,\; 1\; ,\; 1}]} \end{array}\] 
Hence, by Definition 4, we have
\[\overline{X}=[{\rm 0.9982\; ,\; 0.7552\; ,\; 0.7955\; ,\; 0.7456\; ,\; 0.9908\; ,\; 0.9107\; ,\; 1}].\] 
Also, by Definition 3 and Theorem 3, for example, the minimal solutions of $S(a_{1} ,b_{1}^{} )$are obtained as follows:  
\[\underline{X}(1,2)=[{\rm 0\; ,\; 0.7552\; ,\; 0\; ,\; 0\; ,\; 0\; ,\; 0\; ,\; 0}], \underline{X}(1,3)=[{\rm 0\; ,\; 0\; ,\; 0.9341\; ,\; 0\; ,\; 0\; ,\; 0\; ,\; 0}]\] 
Therefore, by Theorem 4, $S(a_{1} ,b_{1} )=[\underline{X}(1,2)\, ,\, \overline{X}\eqref{GrindEQ__1_}\, ]\bigcup [\underline{X}(1,3)\, ,\, \overline{X}\eqref{GrindEQ__1_}\, ]$.

\noindent According to Corollary 2, since $\overline{X}\in S(A,b)$, then the problem is feasible. On the other hand, from Definition 6, we have $\left|\, E\, \right|=24$. Therefore, the number of all vectors $e\in E$ is equal to 24. However, each solution \textbf{$\underline{X}(e)$ }generated by vectors $e\in E$ is not necessary a feasible minimal solution. 

\noindent Additionally, we have $J^{\, \infty } \eqref{GrindEQ__1_}=\{ 1,4,5,6,7\} $, $J^{\, \infty } (2)=\{ 2,3,5,6,7\} $, $J^{\, \infty } (3)=\{ 1,3,4,7\} $, $J^{\, \infty } (4)=\{ 1,2,4,5,6,7\} $ and $J^{\, \infty }(5)=\{ 1,3,5,6,7\} $. So, from the first simplification technique (Lemma 4), ``resetting all the components \textbf{$a_{ij} $ }($i\in I$ and $j\in J^{\, \infty } (i)$) to zeros'' are equivalence operations. Also, by Lemma 5 (second simplification), we can change the value of components \textbf{$a_{13} $}, \textbf{$a_{24} $}, \textbf{$a_{32} $} and \textbf{$a_{52} $} to zeros. For example, since $3\in J(1)\bigcup J(4)$ and ${\rm 0.0404}=b_{4}^{\, p} -b_{1}^{\, p} <w\, (a_{43}^{\, p} -a_{13}^{\, p} )={\rm 0.1183}$, Lemma 5 implies \textbf{$a_{13} =0$}.

\noindent By applying the simplification methods, $\left|\, E\, \right|$ is decreased from 24 to 2.  Therefore, the simplification processes reduced the number of the minimal candidate solutions from 24 to 2, by removing 22 points \textbf{$\underline{X}(e)$}. Indeed, the feasible region has 2 minimal solutions as follows:
\[\begin{array}{l} {e_{1} =[{\rm 2,1,6,3,4}]\Rightarrow \underline{X}(e_{1} )=[{\rm 0.9982\; ,\; 0.7552\; ,\; 0.7955\; ,\; 0.7456\; ,\; 0\; ,\; 0.9107\; ,\; 0}]} \\ {e_{2} =[{\rm 2,1,5,3,4}]\Rightarrow \underline{X}(e_{2} )=[{\rm 0.9982\; ,\; 0.7552\; ,\; 0.7955\; ,\; 0.7456\; ,\; 0.9908\; ,\; 0\; ,\; 0}]} \end{array}\] 
By comparison of the values of the objective function for the minimal solutions, \textbf{$\underline{X}(e_{1} )$ }is optimal for \textbf{$Z_{2} $} (i.e., \textbf{$e^{*} =e_{1} $}). Thus, from Theorem 6, \textbf{$x^{*} =[{\rm 0.9982\; ,\; 0.7552\; ,\; 0.7955\; ,\; 0.7456\; ,\; 0\; ,\; 0.9107\; ,\; 0}]$ }and then \textbf{$Z^{*} =c^{T} x^{*} =-\, {\rm 15.4085}$}.

\noindent 
\vskip 0.2in
\noindent \textbf{Conclusion }
\vskip 0.2in
\noindent In this paper, we proposed an algorithm to solve the linear optimization model constrained with weighted power mean fuzzy relational equalities (WPM-FRE). The feasible solutions set of each WPM-FRE was obtained and their feasibility conditions were described. Based on the foregoing results, the feasible region of the problem is completely resolved. It was shown that the feasible solutions set can be write in terms of a finite number of closed convex cells. Moreover, two simplification operations (depending on the max-WPM composition) were proposed to accelerate the solution of the problem. Finally, a method was introduced for finding the optimal solutions of the problem. 

\noindent \textbf{}

\noindent \textbf{References }
\vskip 0.2in
\noindent [1]. Chang C. W., Shieh B. S., Linear optimization problem constrained by fuzzy max--min relation equations, Information Sciences 234 (2013) 71--79

\noindent [2]. Chen L., Wang P. P., Fuzzy relation equations (i): the general and specialized solving algorithms, Soft Computing 6 (5)(2002) 428-435.

\noindent [3]. Chen L., Wang P. P., Fuzzy relation equations (ii): the branch-point-solutions and the categorized minimal solutions, Soft Computing 11 \eqref{GrindEQ__1_} (2007) 33-40. 

\noindent [4]. Dempe S., Ruziyeva A., On the calculation of a membership function for the solution of a fuzzy linear optimization problem, Fuzzy Sets and Systems 188 (2012) 58-67.

\noindent [5]. Di Nola A., Sessa S., Pedrycz W., Sanchez E., Fuzzy relational Equations and their applications in knowledge engineering, Dordrecht: Kluwer Academic Press, 1989. 

\noindent [6]. Dubey D., Chandra S., Mehra A., Fuzzy linear programming under interval uncertainty based on IFS representation, Fuzzy Sets and Systems 188 (2012) 68-87.

\noindent [7]. Dubois D., Prade H., Fundamentals of Fuzzy Sets, Kluwer, Boston, 2000.

\noindent [8]. Fan Y. R., Huang G. H., Yang A. L., Generalized fuzzy linear programming for decision making under uncertainty: Feasibility of fuzzy solutions and solving approach, Information Sciences 241 (2013) 12-27. 

\noindent [9]. Fang S.C., Li G., Solving fuzzy relational equations with a linear objective function, Fuzzy Sets and Systems 103 (1999) 107-113.

\noindent [10]. Freson S., De Baets B., De Meyer H., Linear optimization with bipolar max--min constraints, Information Sciences 234 (2013) 3--15.

\noindent [11]. Ghodousian A., Optimization of linear problems subjected to the intersection of two fuzzy relational inequalities defined by Dubois-Prade family of t-norms, Information Sciences 503 (2019) 291--306.

\noindent [12]. Ghodousian A., Khorram E., Fuzzy linear optimization in the presence of the fuzzy relation inequality constraints with max-min composition, Information Sciences 178 (2008) 501-519.

\noindent [13]. Ghodousian A., Khorram E., Linear optimization with an arbitrary fuzzy relational inequality, Fuzzy Sets and Systems 206 (2012) 89-102.

\noindent [14]. Ghodousian A., Raeisian Parvari M., A modified PSO algorithm for linear optimization problem subject to the generalized fuzzy relational inequalities with fuzzy constraints (FRI-FC), Information Sciences 418--419 (2017) 317--345.

\noindent [15]. Ghodousian A., Naeeimi M., Babalhavaeji A., Nonlinear optimization problem subjected to fuzzy relational equations defined by Dubois-Prade family of t-norms, Computers \& Industrial Engineering 119 (2018) 167--180.

\noindent [16]. Ghodousian A., Babalhavaeji A., An efficient genetic algorithm for solving nonlinear optimization problems defined with fuzzy relational equations and max-Lukasiewicz composition, Applied Soft Computing 69 (2018) 475--492.

\noindent [17]. Guo F. F., Pang L. P., Meng D., Xia Z. Q., An algorithm for solving optimization problems with fuzzy relational inequality constraints, Information Sciences 252 ( 2013) 20-31.

\noindent [18]. Guo F., Xia Z. Q., An algorithm for solving optimization problems with one linear objective function and finitely many constraints of fuzzy relation inequalities, Fuzzy Optimization and Decision Making 5 (2006) 33-47.

\noindent [19]. Guu S. M., Wu Y. K., Minimizing a linear objective function under a max-t-norm fuzzy relational equation constraint, Fuzzy Sets and Systems 161 (2010) 285-297.

\noindent [20]. Guu S. M., Wu Y. K., Minimizing a linear objective function with fuzzy relation equation constraints, Fuzzy Optimization and Decision Making 12 (2002) 1568-4539.

\noindent [21]. Guu S. M., Wu Y. K., Minimizing an linear objective function under a max-t-norm fuzzy relational equation constraint, Fuzzy Sets and Systems 161 (2010) 285-297.

\noindent [22]. Guu S. M., Wu Y. K., Minimizing a linear objective function with fuzzy relation equation constraints, Fuzzy Optimization and Decision Making 1 (3) (2002) 347-360.

\noindent [23]. Lee H. C., Guu S. M., On the optimal three-tier multimedia streaming services, Fuzzy Optimization and Decision Making 2\eqref{GrindEQ__1_} (2002) 31-39.

\noindent [24]. Li P., Liu Y., Linear optimization with bipolar fuzzy relational equation constraints using lukasiewicz triangular norm, Soft Computing 18 (2014) 1399-1404.

\noindent [25]. Li P., Fang, A survey on fuzzy relational equations, part I: classification and solvability, Fuzzy Optimization and Decision Making 8 (2009) 179-229.

\noindent [26]. Li J. X., Yang S. J., Fuzzy relation inequalities about the data transmission mechanism in bittorrent-like  peer-to-peer file sharing systems, in: Proceedings of the 9${}^{th}$ International Conference on Fuzzy Systems and Knowledge discovery (FSKD 2012), pp. 452-456.  

\noindent [27]. Li P. K., Fang S. C., On the resolution and optimization of a system of fuzzy relational equations with sup-t composition, Fuzzy Optimization and Decision Making 7 (2008) 169-214. 

\noindent [28]. Lin J. L., Wu Y. K., Guu S. M., On fuzzy relational equations and the covering problem, Information Sciences 181 (2011) 2951-2963.

\noindent [29]. Lin J. L., On the relation between fuzzy max-archimedean t-norm relational equations and the covering problem, Fuzzy Sets and Systems 160 (2009) 2328-2344. 

\noindent [30]. Liu C. C., Lur Y. Y., Wu Y. K., Linear optimization of bipolar fuzzy relational equations with max-{\L}ukasiewicz composition, Information Sciences 360 (2016) 149--162.

\noindent [31]. Loetamonphong J., Fang S. C., Optimization of fuzzy relation equations with max-product composition, Fuzzy Sets and Systems 118 (2001) 509-517.

\noindent [32]. Markovskii A. V., On the relation between equations with max-product composition and the covering problem, Fuzzy Sets and Systems 153 (2005) 261-273.

\noindent [33]. Mizumoto M., Zimmermann H. J., Comparison of fuzzy reasoning method, Fuzzy Sets and Systems 8 (1982) 253-283.

\noindent [34]. Peeva K., Resolution of fuzzy relational equations-methods, algorithm and software with applications, Information Sciences 234 (2013) 44-63.

\noindent [35]. Pedrycz W., Granular Computing: Analysis and Design of Intelligent Systems, CRC Press, Boca Raton, 2013.

\noindent [36]. Perfilieva I., Finitary solvability conditions for systems of fuzzy relation equations, Information Sciences 234 (2013)29-43.

\noindent [37]. Qu X. B., Wang X. P., Lei Man-hua. H., Conditions under which the solution sets of fuzzy relational equations  over complete Brouwerian lattices form lattices, Fuzzy Sets and Systems 234 (2014) 34-45.

\noindent [38]. Qu X. B., Wang X. P., Minimization of linear objective functions under the constraints expressed by a system of fuzzy relation equations, Information Sciences 178 (2008) 3482-3490.

\noindent [39]. Sanchez E., Solution in composite fuzzy relation equations: application to medical diagnosis in Brouwerian logic, in: M.M. Gupta. G.N. Saridis, B.R. Games (Eds.), Fuzzy Automata and Decision Processes, North-Holland, New York, 1977, pp. 221-234.

\noindent [40]. Shieh B. S., Infinite fuzzy relation equations with continuous t-norms, Information Sciences 178 (2008) 1961-1967.

\noindent [41]. Shieh B. S., Minimizing a linear objective function under a fuzzy max-t-norm relation equation constraint, Information Sciences 181 (2011) 832-841.

\noindent [42]. Sun F., Wang X. P., Qu X. B., Minimal join decompositions and their applications to fuzzy relation equations over complete Brouwerian lattices, Information Sciences 224 (2013) 143-151. 

\noindent [43]. Sun F., Conditions for the existence of the least solution and minimal solutions to fuzzy relation equations over complete Brouwerian lattices, Information Sciences 205 (2012) 86-92.

\noindent [44]. Wu Y. K., Guu S. M., Minimizing a linear function under a fuzzy max-min relational equation constraints, Fuzzy Sets and Systems 150 (2005) 147-162.

\noindent [45]. Wu Y. K., Guu S. M., An efficient procedure for solving a fuzzy relation equation with max-Archimedean t-norm composition, IEEE Transactions on Fuzzy Systems 16 (2008) 73-84.

\noindent [46]. Wu Y. K., Optimization of fuzzy relational equations with max-av composition, Information Sciences 177 (2007) 4216-4229.

\noindent [47]. Wu Y. K., Guu S. M., Liu J. Y., Reducing the search space of a linear fractional programming problem under fuzzy relational equations with max-Archimedean t-norm composition, Fuzzy Sets and Systems 159 (2008) 3347-3359.

\noindent [48]. Xiong Q. Q., Wang X. P., Fuzzy relational equations on complete Brouwerian lattices, Information Sciences 193 (2012) 141-152.

\noindent [49]. Yang S. J., An algorithm for minimizing a linear objective function subject to the fuzzy relation inequalities with addition-min composition, Fuzzy Sets and Systems 255 (2014) 41-51.

\noindent [50]. Yang X. P., Zhou X. G., Cao B. Y., Latticized linear programming subject to max-product fuzzy relation inequalities with application in wireless communication, Information Sciences\textbf{ }358--359 (2016) 44--55.

\noindent [51]. Yeh C. T., On the minimal solutions of max-min fuzzy relation equations, Fuzzy Sets and Systems 159 (2008) 23-39.  \textbf{}

\end{document}